\documentclass[112pt]{article}
\usepackage{epsf,latexsym}
\usepackage{epsfig}
\usepackage[english]{babel}
\usepackage{amsmath,amssymb,latexsym,amscd}
\usepackage{verbatim}

\title{The Klein bottle in its classical shape:\\
a further step towards a good parametrization}
\author{Gregorio Franzoni}

 \def\R{\mathbb{R}} 
 \def\K{\mathbb{K}} 
 \def\Sph{\mathbb{S}} 

\begin{document}

\maketitle

\noindent \bf Abstract \\
\noindent \rm Together with the M\"obius strip, the Klein bottle
is one of the intriguing objects in the universe of geometry,
sometimes appearing in non-mathematical contexts too. Until now,
several parametrizations of it as a surface immersed in ordinary
three-space have been found, some of which are very elegant and
lead to nice and well understandable shapes. Nevertheless, these
shapes are quite different from the object imagined by F. Klein in
the late $19^{th}$ century: a tube which passes through itself
with the two ends glued together. Parametrizations for this
version of the Klein bottle do exist, but they are not fully
satisfactory for some reasons. We discuss some of the existing
representations and propose two new immersions of the Klein bottle
in $\R^3$, which are intended to be a step towards a canonical
expression of this surface in the shape imagined by its first
discoverer.

\section{Introduction}

The Klein bottle is a topological object that can be defined as
$\K=\R^2 \big/ \sim$, where $\sim$ is the equivalence relation
\mbox{$(u,v)\sim (u',v')\iff (u',v')=(u+2k\pi,(-1)^k v+ 2h \pi)$}.
It is a well known fact that $\K$ is a genus 2 non-orientable
closed surface and that it is non-embeddable in $\R^3$. It is
possible, however, to immerse it in $\R^3$, that is, to map it in
$\R^3$ obtaining an image with no singular points. To give an
immersion of $\K$ in $\R^3$ it suffices to define, on the
fundamental square $[0, 2 \pi] \times [0, 2 \pi]$, an immersion
that passes to the quotient with respect to the relation $\sim$.
Among the existing known immersions of $\K$ in $\R^3$, we recall
the two following:
\begin{equation} \label{lemnisc_param}
\begin{array}{ll}
\mathrm{Kb1}(u,v):\left\{
\begin{array}{ll}
x=(a+\cos(\frac{u}{2})\sin v-\sin(\frac{u}{2})\sin(2v) \cos u\\
y=(a+\cos(\frac{u}{2}) \sin v-\sin(\frac{u}{2})\sin(2v)) \sin u\\
z=\sin(\frac{u}{2}) \sin v+\cos(\frac{u}{2})\sin(2v)
\end{array} \right.\\[8mm]
\end{array}
\end{equation}
\begin{figure}[!h]
  \begin{center}
   \epsfig{file=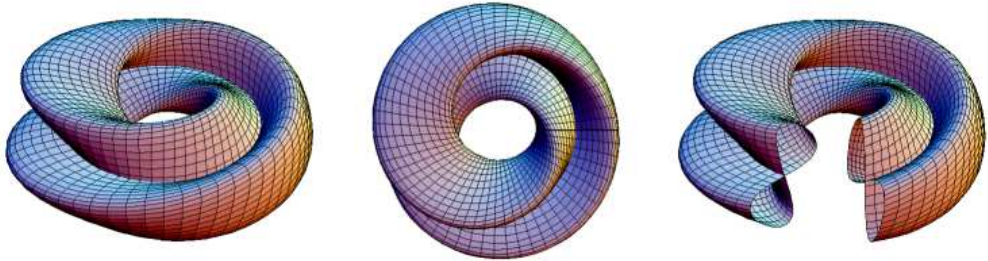, width=\textwidth}
    \caption{Three views of a Klein bottle generated by moving a \emph{lemniscate}.}
   \label{figure8_kleinb}
\end{center}
\end{figure}
\begin{equation} \label{sieben_param}
\begin{array}{ll}
\mathrm{Kb2}(u,v):\left\{
\begin{array}{ll}
x=\cos 2u \ \sin v /(1-(\sin u \ \cos u+\sin 2u \ \sin v)/\sqrt 2)\\
y=(\sin 2u \ \sin v - \sin u \ \cos v)/\sqrt 2
(1-(\sin u \ \cos u+\sin 2u \ \sin v)/\sqrt 2)\\
z=\cos u \ \cos v /(1-(\sin u \ \cos u+\sin 2u \ \sin v)/\sqrt 2)
\end{array}  \right.\\[8mm]
\end{array}
\end{equation}
\begin{figure}[!h]
  \begin{center}
   \epsfig{file=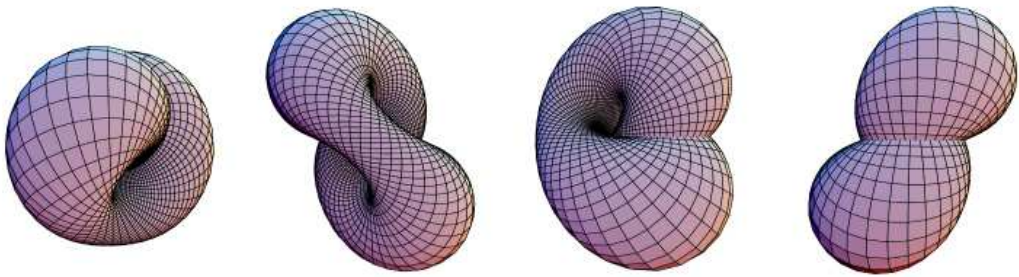, width=\textwidth}
    \caption{Klein bottle as a one-parameter family of circles,
    four views from coordinate directions.}
   \label{paramlawson_multi}
\end{center}
\end{figure}
\begin{figure}[!h]
  \begin{center}
   \epsfig{file=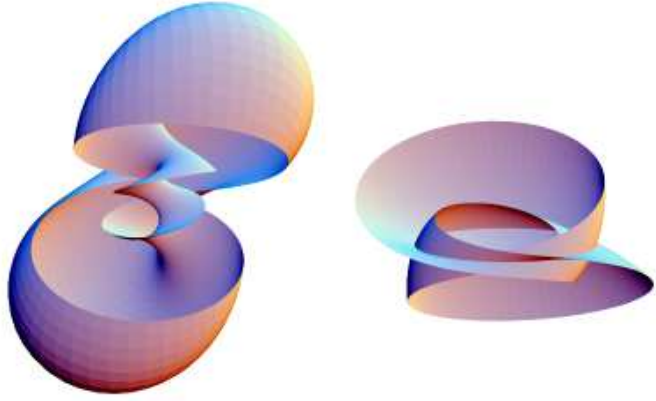, width=0.65\textwidth}
    \caption{Klein bottle as a one-parameter family of circles,
    interior view.}
   \label{cutlawson_multi}
\end{center}
\end{figure}
The first one can be obtained by moving a \emph{lemniscate} around
a circle while it rotates, in its plane, half a turn around its
center (see, for example, \cite{ap}, \cite{ban}, \cite{pin}). The
second one can be seen as a one-parameter family of circles in
space and its image is also the zero set of a polynomial of degree
eight (see \cite{ap}). Moreover, in the work of H. B. Lawson
\cite{law}, it arises as particular case as the image, by
stereographic projection, of a minimal surface on $\Sph^3$. The
two versions of the Klein bottle described above could be
considered \emph{canonical} for the simplicity of their
geometrical construction and of their formulas, and for the
clearness of the resulting shapes.

But if we think of the bottle as it was conceived by F. Klein in
1882 (\cite[\S 23]{kle}, see also \cite[pp. 308--310]{hilb}):
\emph{tucking a rubber hose, making it to penetrate to itself, and
then smoothly gluing the two ends together}, they are far from
such a visual idea. On the other hand, known formulas that lead to
the bottle in its original shape are quite complicated and don't
have the elegance of \eqref{figure8_kleinb} and
\eqref{sieben_param}. In the following we will recall two
definitions of the Klein bottle in its first shape, then we will
propose a further step towards what we would like to call a
canonical parametrization for it as a self-penetrating tube.

\section{Two classical looking definitions of the Klein bottle}
In this paragraph we recall two parametrizations of the bottle due
to S. Dickson \cite{dic} and to M.~Trott \cite{tro} respectively.
Both surfaces are realized as tubes around plane curves, according
to the following scheme: let $\alpha(t)=(x(t),y(t))$, $t \in
[a,b]$ be a curve lying on the $xy$ plane satisfying
$\|\alpha^{\prime}(t)\| \neq 0$, $\mathbf{k}=(0,0,1)$ the $z$-axis
unit vector,
$\mathbf{T}=\frac{\alpha^{\prime}}{\|\alpha^{\prime}\|}$ the unit
tangent vector field of $\alpha(t)$ and $J$ the canonical complex
structure of the $xy$ plane, i.\hspace{1mm}e. the linear map from
$\R^2$ to itself defined by $J(v_1,v_2)=(-v_2,v_1)$. Then the
couple of unit vectors $(J(\mathbf{T}),\mathbf{k})$ is always
orthogonal to $\alpha^{\prime}(t)$ along $\alpha(t)$ and can be
used to construct a tube around $\alpha(t)$ as follows:

\begin{equation} \label{scheme_tubecurve}
\begin{array}{ll}
\mathrm{Tube}(t,\theta)=\alpha(t)+r(t)\big(\cos \theta \
J(\mathbf{T})+
 \sin \theta \ \mathbf{k} \big)\\[2mm]
(t,\theta)\in [a,b]\times[0,2\pi]
\end{array}
\end{equation}

\noindent where the scalar continuous function $r(t)$ represents
the radius of the tube.

\begin{figure}[!h]
  \begin{center}
   \epsfig{file=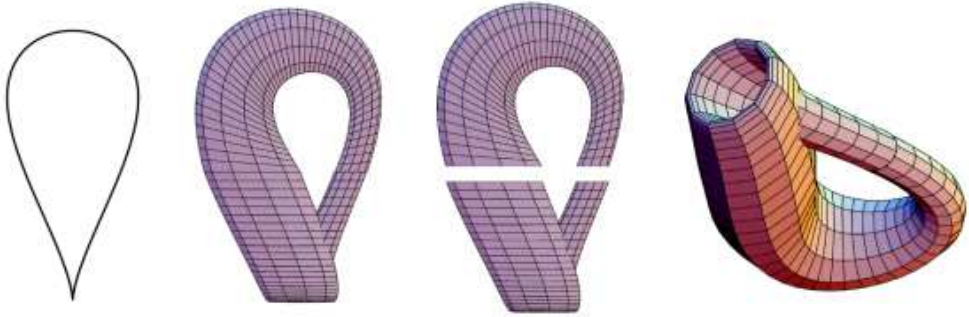, width=\textwidth}
    \caption{Klein bottle according to S. Dickson's definition;
    on the left, the central curve.}
   \label{dicksonbottle}
\end{center}
\end{figure}

Dickson's bottle is built up by that scheme, except for the choice
of the moving pair of vector fields, which is not orthogonal to
the curve. Here is its parametrization
\begin{equation} \label{dicks_param} \mathrm{Kb3}(u,v)=\left\{
\begin{array}{ll}
x= \left\{
\begin{array}{ll}
6\cos u (1+\sin u)+4(1-\frac{1}{2}\cos u)\cos u \cos v
\qquad \textrm{for} \quad 0 \leq u \leq \pi\\[3mm]
6\cos u (1+\sin u)+4(1-\frac{1}{2}\cos u)\cos(v+\pi)\qquad
\textrm{for} \quad \pi<u \leq 2\pi
\end{array} \right.\\[6mm]
y= \left\{
\begin{array}{ll}
16 \sin u + 4(1-\frac{1}{2}\cos u)\sin u \cos v \qquad
\textrm{for} \quad 0 \leq u \leq \pi\\[3mm]
16 \sin u \qquad \textrm{for} \quad \pi<u \leq 2\pi
\end{array} \right.\\[6mm]
z=4(1-\frac{1}{2}\cos u)\sin v
\end{array} \right.
\end{equation}
Parametrization \eqref{dicks_param} defines two distinct tubes,
the first one being built up on a frame which moves along the
central curve remaining parallel to the $xz$ plane, the second one
connecting the two ends of the first tube through a rotation of
$\pi$ of the moving frame. The union of these two parts turns out
to be a well-looking object (Fig. \ref{dicksonbottle}), which
renders properly the idea outlined by F. Klein.

\begin{figure}[!h]
  \begin{center}
   \epsfig{file=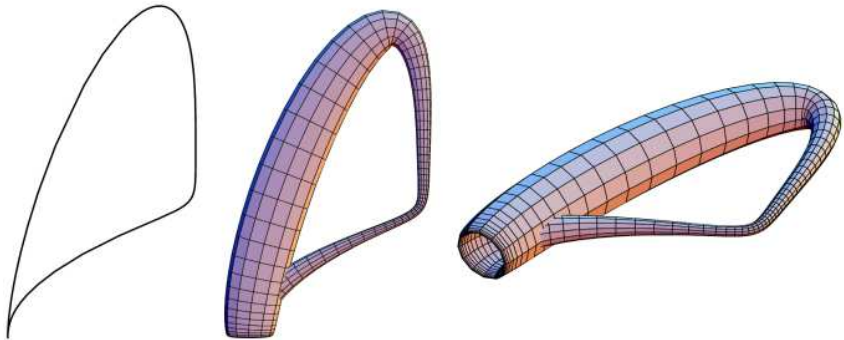, width=0.8\textwidth}
    \caption{Klein bottle according to M. Trott's definition;
    on the left, the central curve.}
   \label{trottbottle}
\end{center}
\end{figure}

\noindent Note that the curve implemented in this construction
\begin{equation}
\alpha(t): \left\{
\begin{array}{ll}
x = 6 \cos t (1+\sin t)\\
y = 16 \sin t
\end{array} \right.
\end{equation}
is a \emph{piriform}, a well known curve (see for example
\cite{ca-gr} or \cite{cem}) whose parametrization is
\begin{equation}
\textrm{piriform}(t): \left\{
\begin{array}{ll}
x = a (1+\sin t)\\
y = b \cos t (1+\sin t)
\end{array} \right.\\[4mm]
\end{equation}
It can be easily proven, and also guessed by looking at Figure
\ref{dicksonbottle}, that \eqref{dicks_param} does not define an
immersion, as the two tubes do not meet tangentwise along the
common boundaries. And, of course, it would be better to get the
\emph{whole} surface as the image of a single parametrization,
with no use of inequalities as in \eqref{dicks_param}. This has
been done by M. Trott, who follows closely the scheme defined by
\eqref{scheme_tubecurve} and puts some constraints on the central
curve (that he calls \emph{directrix}) and on the radius, that can
be resumed in the following:
\begin{equation} \label{constraints}
\begin{array}{ll}
i)& \ \alpha(a)=\alpha(b)\\[2mm]
ii)& \ \alpha'(a)=-\alpha'(b)\\[2mm]
iii)& \ r(a)=r(b)\\[2mm]
iv)& \ r'(a)=r'(b)= \pm \infty\\
\end{array}
\end{equation}
Conditions $i)$, $ii)$ and $iii)$ mean that the two ends of the
tube must be coincident, while $iv)$ means that they must meet
tangentwise. The curve and the radius he uses are
\begin{equation} \label{trott_param}
\begin{array}{ll}
\displaystyle \beta(t)=\left( \frac{1}{t^4+1} \, , \,
\frac{t^2+t+1}{t^4+1}\right)  \hspace{1cm} t \in (-\infty, +\infty)\\[6mm]
\displaystyle r(t)= \frac{84 t^4+ 56 t^3 + 21 t^2 + 21 t + 24}{672
(1+t^4)}
\end{array}
\end{equation}
and the resulting image is shown in Figure \ref{trottbottle}.
Equations \eqref{trott_param} define an immersion, but the
resulting shape is somehow edgy, because of the choice of a
directrix whose curvature has a non-smooth behavior. Moreover, as
$t$ ranges on an open interval, when one tries to plot the surface
there is a missing strip in correspondence of a neighborhood of
the cusp (Trott uses $t \in [-20, +20]$ in his plots).
\section{A new description}
\begin{figure}[!h]
  \begin{center}
   \epsfig{file=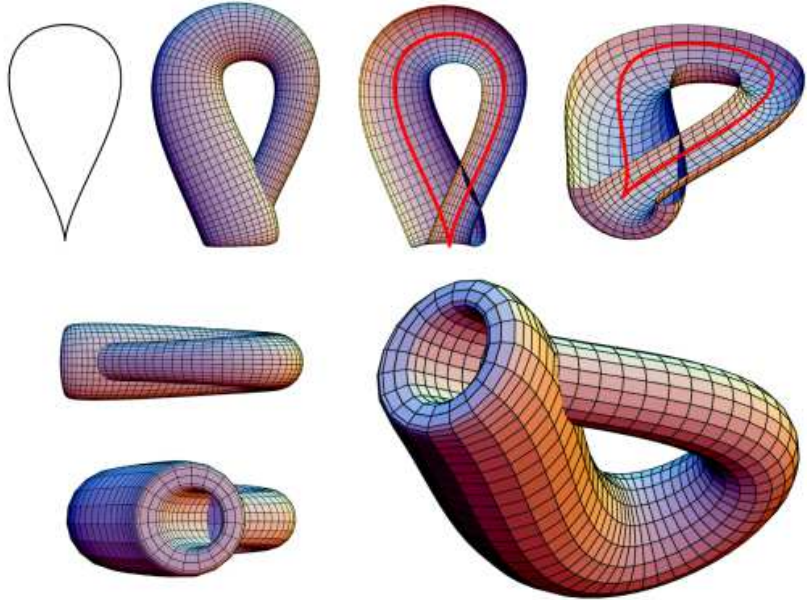, width=0.80\textwidth}
    \caption{The piriform curve and a tube around it: an
    immersion of the Klein bottle in $\R^3$.}
    \label{franzonikb}
   \end{center}
\end{figure}

Starting from the two constructions described in the previous
paragraph, it is natural trying to construct a new surface by
taking the best features from both: the beautiful and symmetric
directrix of the first one and the rigorous geometric scheme of
the second one. In order to use the piriform as a directrix for
our tube, we re-parametrize it to make it start and end at the
cusp:
\begin{equation} \label{pirif_reparam}
\begin{array}{ll}
\gamma(t)=\left\{
\begin{array}{ll}
a(1-\cos t)\\
b \sin t (1-\cos t)
\end{array} \right.\\[4mm]
\hspace{1.5cm} t \in (0,2 \pi)
\end{array}
\end{equation}
A suitable radius, which satisfies $iii)$ and $iv)$ of
\eqref{constraints}, is for example
\begin{equation} \label{myradius} r(t) =
c-d(t-\pi)\sqrt{t(2\pi-t)}
\end{equation}
Parameters $c$ and $d$ affect respectively the radius of the whole
tube and the difference between its minimum and maximum value. The
resulting plot, with $(a,b,c,d)=(20,8,\frac{11}{2},\frac{2}{5})$
and $(t, \theta) \in (0,2 \pi) \times [0,2 \pi]$ is shown in
Figure \ref{franzonikb}.
\section{Some remarks}
Although, in our opinion, the described result is rather
satisfactory, there are some facts we want to point out. First,
the parametrization of our surface, extensively written, has a
long and complicated expression. Secondly, similarly to what
happens with Trott's parametrization, the image of the immersion
fails to be closed because it misses a circle at the cusp of the
directrix, as $\|{\gamma}^{\prime}\|$ vanishes at $t=0$ and
$t=\pi$, while the used scheme needs $\|\gamma'\|$ to be nonzero
everywhere. A way to get rid of this issue is using one half of
the Dumbbell curve (see \cite{cem}) as a directrix. This is a
famous sextic curve which also has the following parametrization:
\begin{equation}
\textrm{dumbbell}(t): \left\{
\begin{array}{ll}
x = \sin t\\
y = \sin ^2 t \cos t
\end{array} \right.
\hspace{6mm} t \in [0,2 \pi]
\end{equation}
\begin{figure}[!h]
  \begin{center}
   \epsfig{file=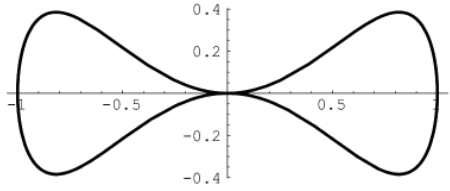, width=0.45\textwidth}
    \caption{Dumbbell curve.}
   \label{dumbbellgraphics}
\end{center}
\end{figure}
If $t$ ranges in $I=[0, \pi]$, one obtains a curve that satisfies
the first three conditions of \eqref{constraints} and whose
tangent vector is well-defined for all $t \in I$, so it is
possible to use a closed rectangle as a domain for the immersion,
obtaining a closed image. By using a stretched Dumbbell curve
$\alpha(t)=(5 \sin t, 2 \sin^2 t \cos t)$ as directrix and $r(t) =
\frac{1}{2}-\frac{1}{30} (2t-\pi)\sqrt{2t(2\pi-2t)}$ as radius
function, with $t \in [0, \pi]$, we obtain another example of
immersion of the Klein bottle, more suitable to be plotted with a
computer (see Figure \ref{kleinb_dumbbell}).
\begin{figure}[!h]
  \begin{center}
   \epsfig{file=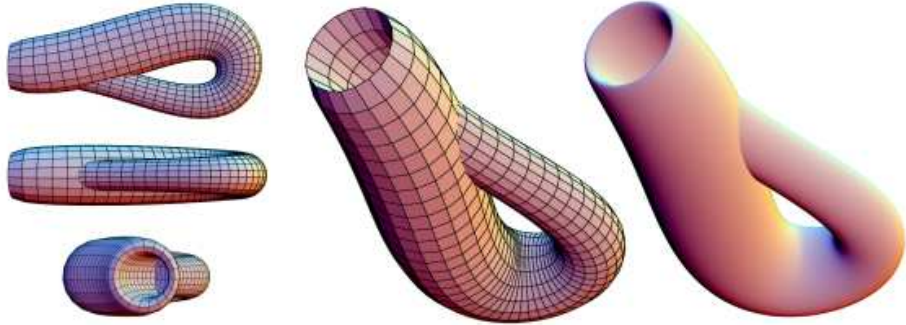, width=0.9\textwidth}
    \caption{Some views of the Klein bottle as a tube around
    Dumbbell curve. On the right, a high definition plot.}
   \label{kleinb_dumbbell}
\end{center}
\end{figure}
\section{Conclusion}
We defined two immersions of the Klein bottle in $\R^3$ in the
shape outlined by F. Klein in 1882, with a reasonably good
appearance. The mathematical expression of both is still too
complicated and far from the elegance of versions like
\eqref{lemnisc_param} and \eqref{sieben_param}. The present work
is intended to be a further step towards an immersion of this
famous surface in $\R^3$ which we would like to call ``canonical''
from both a mathematical and an aesthetic point of view.

\end{document}